\documentclass[11pt,twoside]{article}

\usepackage{amsmath}
\usepackage{amsfonts}
\usepackage{amssymb}
\newcommand{\RR}{\mathbb R}

\newcommand{\email}[1]{{\small E-mail: {\textsf {#1}}}}
\newtheorem{defi}{Definition}

\newtheorem{theo}[defi]{Theorem}
\newtheorem{Le}[defi]{Lemma}

\newtheorem{Rem}[defi]{Remark}

\newcommand\qed{\hfill $\square$}

\begin{document}

\title{Finite time blow-up for radially symmetric solutions to a critical quasilinear Smoluchowski-Poisson system} 

\author{Tomasz Cie\'slak\footnote{Institute of Applied Mathematics and Mechanics, University of Warsaw, Banacha 2, 02-097 Warszawa, Poland. \email{T.Cieslak@impan.gov.pl}} \kern8pt \& \kern8pt
Philippe Lauren\c cot\footnote{Institut de 
Math\'ematiques de Toulouse, CNRS UMR~5219, Universit\'e de Toulouse, 118 route de Narbonne, F--31062 Toulouse Cedex 9, 
France. \email{Philippe.Laurencot@math.univ-toulouse.fr}}}
\date{\today}
\maketitle

\begin{abstract}
Finite time blow-up is shown to occur for radially symmetric solutions to a critical quasilinear Smoluchowski-Poisson system provided that the mass of the initial condition exceeds an explicit threshold. In the supercritical case, blow-up is shown to take place for any positive mass. The proof relies on a novel identity of virial type.
\end{abstract}

\section{Introduction}\label{i}


We study the occurrence of blow-up in finite time for radially symmetric solutions to a generalized Smoluchowski-Poisson system
\begin{eqnarray}
\label{he1}
\partial_t{u} &=& {\rm div }\left( a(u)\ \nabla{u}-u\ \nabla{v} \right) \;\;\mbox{in}\;\;(0,\infty)\times B(0,1),\\
\label{he12}
0&=&\Delta v + u - M \;\;\mbox{in}\;\;(0,\infty)\times B(0,1),\\
\label{he2}
\partial_\nu u & = & \partial_\nu v = 0\;\;\mbox{on} \;\; (0,\infty)\times \partial B(0,1),\\
\label{he3}
u(0)&=& u_0\ge 0 \;\;\mbox{in} \;\;B(0,1),\;\;\int_{B(0,1)} v(t,x)\ dx=0\;\;\mbox{for any}\;\; t\in (0,\infty),
\end{eqnarray}
where $B(0,1)$ denotes the unit ball of $\RR^n$, $n\ge 2$, and $M$ the mean value of $u_0$. The diffusion coefficient $a$ belongs to $\mathcal{C}^2([0,\infty))$ and is assumed to be positive for simplicity (see Remark~\ref{rezz2} below). The system \eqref{he1}-\eqref{he3} arises in astrophysics as a model of self-gravitating Langevin particles \cite{Ch03} and in biology \cite{KS70} where it is also known as the parabolic-elliptic Keller-Segel chemotaxis model.

A fundamental property of solutions to \eqref{he1}-\eqref{he3} is that 
\begin{equation}
\label{zz1}
u(t)\ge 0\,, \quad \int_{B(0,1)} u(t,x)\ dx = M\ |B(0,1)|\,, \;\;\mbox{ and }\;\; \int_{B(0,1)} v(t,x)\ dx = 0
\end{equation}
for $t>0$, which readily follows from \eqref{he1}, \eqref{he2}, the comparison principle, the non-negativity of $u_0$, and the definition of $v$. It is by now well-known that, if $a\equiv 1$, there are radially symmetric initial data $u_0$ for which the corresponding solution to \eqref{he1}-\eqref{he3} blows up in finite time and this singular phenomenon may happen for any $M>0$ if $n\ge 3$ but only for $M>8\pi$ if $n=2$ \cite{JL92,Nagai}. More recently, it was shown in \cite{cw_1} that there is a critical exponent for the nonlinear diffusion coefficient $a$ which separates two different behaviours for the solutions to \eqref{he1}-\eqref{he3}: on the one hand, if $a(z)\ge c\ (1+z)^ \alpha $ and $\alpha >(n-2)/n$, there is a unique global classical solution to \eqref{he1}-\eqref{he3} for any non-negative initial condition $u_0\in L^\infty(B(0,1))$ (and this is actually true for a general smooth bounded domain of $\RR^n$, $n\ge 1$). On the other hand, if $a(z)\le c\ (1+z)^ \alpha $ and $\alpha <(n-2)/n$, radially symmetric solutions to \eqref{he1}-\eqref{he3} blowing up in finite time are constructed in \cite{cw_1}. 

Except for $n=2$, the critical case $\alpha =(n-2)/n$ is not covered in \cite{cw_1} and the purpose of this note is to fill this gap: indeed, the main outcome of our analysis is that, if $a(z)\le c\ (1+z)^{(n-2)/n}$, there are solutions to \eqref{he1}-\eqref{he3} blowing up in finite time when $M$ exceeds an explicit threshold. As a by-product of our analysis, we also establish an alternative and simpler proof of the blow-up result in \cite{cw_1} for $\alpha\in[0,(n-2)/n)$. Indeed, the construction of solutions to \eqref{he1}-\eqref{he3} blowing up in finite time performed in \cite{cw_1} relies on the possibility of reducing \eqref{he1}-\eqref{he3} to a single parabolic equation. The approach used in this paper is completely different and relies on the derivation of a differential inequality of virial type which cannot hold true for all times. When $a\equiv 1$, this technique is used in \cite{Nagai} where it is shown that the moment of order $n$ of $u$ satisfies a differential inequality which contradicts the non-negativity of $u$ after a finite time. Seemingly, the moment of order $n$ of $u$ does not give valuable information when the diffusion is nonlinear and we introduce nonlinear functions of $u$ to be able to handle this case.  We finally point out that the above results are only valid for $n\ge 2$: the situation is qualitatively different in the one dimensional case $n=1$ and will be considered in a separate paper \cite{cl1}. 

\section{Finite time blow-up}\label{ftbu}


We first introduce some notation: for $(\alpha,p)\in [0,\infty)\times (1,\infty)$, we define
\begin{equation}
\label{zz2}
\kappa_p(\alpha) := \frac{(p-1)}{(\alpha +1)(p+ \alpha)}\ \left( \frac{2(n-1)}{p-1} \right)^{\alpha +1}\ (np)^{(n-2-\alpha n)/n}\,.
\end{equation}
We also define $\bar{m}_p(f)$ for $p\ge 1$ and $f\in L^\infty(0,1)$ by
\begin{equation}
\label{zz3}
\bar{m}_p(f) := \frac{1}{p}\ \int_0^1 \left( \int_r^1 f(\rho)\ \rho^{n-1}\ d\rho \right)^p\ r^{n-1}\ dr\,.
\end{equation}

Our main result then reads as follows.

\begin{theo}\label{main}
Assume that there are $\alpha\in [0,(n-2)/n]$, and positive real numbers $c_1> 0$ and $c_2>0$ such that
\begin{equation}
\label{zz4}
0<a(z)\le c_1\ z^\alpha + c_2 \;\;\mbox{ for }\;\; z\ge 0\,.
\end{equation}
Let $M>0$ and consider a non-negative radially symmetric function $u_0\in L^\infty(B(0,1))$ such that $\|u_0\|_1=M$. Assume further that $E_{M,p}(\bar{m}_p(u_0))<0$ for some $p>1$, where
\begin{eqnarray}
\nonumber
E_{M,p}(z) & := & c_1\ \kappa_p(\alpha)\ \left( \frac{M}{n} \right)^{(2p+n\alpha(p+1))/n}\ z^{(n-2-\alpha n)/n} + c_2\ \kappa_p(0)\ \left( \frac{M}{n} \right)^{(2p)/n}\ z^{(n-2)/n} \\
\label{zz5}
& + & M\ z - \frac{1}{p(p+1)}\ \left( \frac{M}{n} \right)^{p+1}\,, \quad z\ge 0\,.
\end{eqnarray}
Then the system \eqref{he1}-\eqref{he3} has a unique maximal classical solution $(u,v)$ with finite maximal existence time $T_{max}\in (0,\infty)$ and $\|u(t)\|_\infty\longrightarrow \infty$ as $t\rightarrow T_{max}$.
\end{theo}

\begin{Rem}\label{Rem}
There are initial data $u_0$ for which $E_{M,p}(\bar{m}_p(u_0))<0$ for some $p>1$. Indeed, observe that $E_{M,p}(0)<0$ for all $M>0$ if $\alpha\in [0,(n-2)/n)$ and for $M>n \left( c_1 p(p+1) \kappa_p((n-2)/n) \right)^{n/2}$ if $\alpha=(n-2)/n$. It is then sufficient to take $u_0$ sufficiently concentrated near $x=0$ so that $\bar{m}_p(u_0)$ is close to zero (for instance, $u_0=M \delta^{-n}\ \mathbf{1}_{B(0,\delta)}$ for $\delta>0$ sufficiently small).
\end{Rem}

\noindent\textbf{Proof of Theorem~\ref{main}.}\\
By \cite[Theorem 1.3]{cw_1}, there exist a maximal existence time $T_{max}\in (0,\infty]$ and a unique radially symmetric classical solution $(u,v)\in \mathcal{C}([0,T_{max});L^2(B(0,1);\RR^2))\cap \mathcal{C}^{1,2}((0,T_{max})\times B(0,1);\RR^2)$ to \eqref{he1}-\eqref{he3} satisfying \eqref{zz1} for $t\in [0,T_{max})$. Moreover, if $T_{max}<\infty$ then $\left\|u(t)\right\|_\infty \longrightarrow \infty$ as $t\rightarrow T_{max}$.

We introduce 
$$
U(t,r) := \frac{1}{n |B(0,1)|}\ \int_{B(0,r)} u(t,x)\ dx \;\;\mbox{ and }\;\; m_{p}(t):=\frac{1}{p}\ \int_0^1\left(\frac{M}{n}-U(t,r)\right)^p\ r^{n-1}\ dr
$$
for $(t,r)\in [0,T_{max})\times [0,1]$ and derive the following identity of virial type. 

\begin{Le}\label{lezz6} Let $A$ be defined by $A'=a$ and $A(0)=0$. Then 
\begin{equation}
\label{zz7}
\frac{dm_p}{dt} = M\ m_p - \frac{1}{p(p+1)}\ \left( \frac{M}{n} \right)^{p+1} + \mathcal{R}_p(u)
\end{equation}
with 
$$
\mathcal{R}_p(u) :=  \int_0^1 r^{2n-3}\ \left( \frac{M}{n} - U \right)^{p-2}\  \left( 2(n-1)\ \left( \frac{M}{n} - U \right) - (p-1)\ r^n\ u \right)\ A(u)\ dr\,.
$$
\end{Le}

\noindent\textbf{Proof.} Integrating \eqref{he1} gives that 
$U$ solves 
\[
\partial_t U=r^{n-1}\ \partial_r A(u) + u\ U - \frac{M}{n}\ r^n\ u \;\;\mbox{ with }\;\; U(t,0)=\frac{M}{n}-U(t,1)=0\,.
\]
Consequently,
\begin{eqnarray*}
\frac{d m_p}{dt} & = & \int_0^1 \left( \frac{M}{n}-U \right)^{p-1}\ \left[ \frac{M}{n}\ (r^n -1)\ \partial_r U + \left( \frac{M}{n} - U \right)\ \partial_r U - r^{2(n-1)}\ \partial_r A(u) \right]\ dr \\
& = & -\frac{1}{p}\ \left( \frac{M}{n} \right)^{p+1} + M\ m_p + \frac{1}{p+1}\ \left( \frac{M}{n} \right)^{p+1} \\
& + & \int_0^1 \left[ 2(n-1)\ r^{2n-3}\ \left( \frac{M}{n}-U\right)^{p-1}\ - (p-1)\ r^{2(n-1)}\ \left( \frac{M}{n}-U \right)^{p-2}\ \partial_r U \right]\ A(u)\ dr\,,
\end{eqnarray*}
and hence \eqref{zz7}. \qed

\medskip

The next step is to estimate $\mathcal{R}_p(u)$ in terms of $m_p$. To this end, we notice that the assumption \eqref{zz4} on $a$ warrants that 
\begin{equation}\label{dr}
0 \leq A(z)\leq \frac{c_1}{1+ \alpha}\ z^{1+ \alpha} + c_2\ z\,, \quad z\ge 0\,.
\end{equation}
In view of \eqref{dr}, we have
\begin{eqnarray*}
& & \left[ 2(n-1)\ \left(\frac{M}{n}-U\right)-(p-1)\ r^n\ u\right]\ A(u) \\
& \le & \max{\left\{ 2(n-1)\ \left( \frac{M}{n}-U \right)-(p-1)\ r^n\ u , 0 \right\}}\ \left( \frac{c_1}{1+ \alpha}\ u^{1+ \alpha} + c_2\ u \right) \\
& \le & \max{\left\{ 2(n-1)\ \left( \frac{M}{n}-U \right)-(p-1)\ r^n\ u , 0 \right\}}\  \left( \frac{c_1}{1+ \alpha}\ \left[ \frac{2 (n-1)}{(p-1) r^n}\ \left( \frac{M}{n}-U \right) \right]^ \alpha + c_2 \right)\ u \\
& \le & \frac{c_1 (p-1)}{1+ \alpha}\ \left( \frac{2 (n-1)}{p-1} \right)^{1+ \alpha}\ \left( \frac{M}{n}-U \right)^{1+ \alpha}\ r^{-n \alpha}\ u + 2(n-1) c_2\ \left( \frac{M}{n}-U \right)\ u\,,
\end{eqnarray*}
and thus
\begin{eqnarray}
\nonumber
\mathcal{R}_p(u) & \le & \frac{c_1 (p-1)}{1+ \alpha}\ \left( \frac{2 (n-1)}{p-1} \right)^{1+ \alpha}\ \int_0^1 r^{n-2-\alpha n}\ \left( \frac{M}{n}-U \right)^{p+ \alpha-1}\ \partial_r U\ dr \\
\label{zz8}
& + & 2(n-1) c_2\ \int_0^1 r^{n-2}\ \left( \frac{M}{n}-U \right)^{p-1}\ \partial_r U\ dr\,.
\end{eqnarray}
Since $\alpha\in [0,(n-2)/n]$, the function $r\mapsto r^{(n-2-\alpha n)/n}$ is  concave and we infer from the Jensen inequality (with measure $((M/n)-U)^{p+ \alpha-1}\ \partial_r U\ dr$) that
\begin{eqnarray*}
& & \int_0^1 r^{n-2-\alpha n}\ \left( \frac{M}{n}-U \right)^{p+ \alpha-1}\ \partial_r U\ dr \\
& \le & \left( \frac{1}{p+ \alpha}\ \left( \frac{M}{n} \right)^{p+ \alpha} \right)^{(2+ \alpha n)/n}\ \left( \int_0^1 r^n\ \left( \frac{M}{n}-U \right)^{p+ \alpha-1}\ \partial_r U\ dr \right)^{(n-2-\alpha n)/n}\\
& \le & \left( \frac{1}{p+ \alpha}\ \left( \frac{M}{n} \right)^{p+ \alpha} \right)^{(2+ \alpha n)/n}\ \left( \frac{n}{p+ \alpha}\ \int_0^1 r^{n-1}\ \left( \frac{M}{n}-U \right)^{p+ \alpha}\ \ dr \right)^{(n-2-\alpha n)/n}\\
& \le & \frac{(np)^{(n-2-\alpha n)/n}}{p+ \alpha}\ \left( \frac{M}{n} \right)^{(2p+ \alpha n(p+1))/n}\ m_p^{(n-2-\alpha n)/n}\,.
\end{eqnarray*}
Arguing in a similar way to estimate the second integral in the right-hand side of \eqref{zz8}, we deduce from \eqref{zz8} that 
$$
\mathcal{R}_p(u) \le E_{M,p}(m_p) - M\ m_p + \frac{1}{p(p+1)}\ \left( \frac{M}{n} \right)^{p+1}\,.
$$
Inserting this estimate in \eqref{zz7} we arrive at
\begin{equation}\label{czw}
\frac{d m_p(t)}{dt} \leq E_{M,p}(m_p(t)) \;\;\mbox{ for }\;\; t\in [0,T_{max})\,.
\end{equation}

Assume now for contradiction that $T_{max}=\infty$. Since $z\mapsto E_{M,p}(z)$ is an increasing function and $m_p(0)=\bar{m}_p(u_0)$, we realize that, as soon as  $E_{M,p}(m_p(0))<0$, we have $m_p(t_0)=0$ for some $t_0\in (0,\infty)$. Thus, $U(t_0,r)=M/n$ for all $r\in [0,1]$ which contradicts the fact that $U(t_0,0)=0$.  Consequently, $T_{max}<\infty$ and the proof is complete. \qed

\begin{Rem}\label{rezz1}
If $a\equiv 1$ (i.e. $\alpha =0$), Lemma~\ref{lezz6} is also valid for $p=1$ and $m_1$ coincides with the moment of order $n$ used in \cite{Nagai}. Not surprisingly, if $n=2$ and $a\equiv 1$, we have $E_{M,p}(0)<0$ if $M>4(p+1)$ which converges to $8$ as $p$ approaches $1$ and we recover the well-known threshold condition $\|u_0\|_1>8\pi$ for finite time blow-up to occur in the parabolic-elliptic Keller-Segel system \cite{Nagai}.
\end{Rem}

\begin{Rem}\label{rezz2}
The requirement $a>0$ is only used to have classical solutions to \eqref{he1}-\eqref{he3} but does not play any role in the blow-up condition and the identity of virial type (Lemma~\ref{lezz6}). Thus, Theorem~\ref{main} remains valid if the diffusion is degenerate (for instance, $a(r)=m r^{m-1}$ with $m\in [1,2(n-1)/n]$) provided an appropriate notion of weak solutions is available. 
\end{Rem}

\textbf{Acknowledgement.} This paper was prepared during T.~Cie\'slak's one-month visit at the Institut de Math\'ematiques de Toulouse, Universit\'e Paul Sabatier. T. Cie\'slak would like to express his gratitude for the invitation, support, and hospitality.


\begin{thebibliography}{20}

\bibitem{Ch03}
P.-H. Chavanis, \textit{Generalized thermodynamics and Fokker-Planck equations. Applications to stellar dynamics and two-dimensional turbulence}, Phys. Rev. E \textbf{68} (2003), 036108.

\bibitem{cl1}
T. Cie\'slak and Ph. Lauren\c{c}ot, in preparation.

\bibitem{cw_1}
T. Cie\'slak and M. Winkler, \textit{Finite-time blow-up in a quasilinear system of chemotaxis}, Nonlinearity \textbf{21} (2008), 1057--1076.

\bibitem{JL92}
W.~J\"{a}ger and S.~Luckhaus, \textit{On explosions of
solutions to a system of partial differential equations modelling
chemotaxis}, Trans. Amer. Math. Soc. \textbf{329} (1992), 819--824.

\bibitem{KS70}
E.F. Keller and L.A. Segel, \textit{Initiation of slide mold aggregation viewed as an instability}, J. Theor. Biol. \textbf{26} (1970), 399--415.

\bibitem{Nagai}
T. Nagai, \textit{Blow-up of radially symmetric solutions to a chemotaxis system}, Adv. Math. Sci. Appl. {\bf 5} (1995), 581--601.

\end{thebibliography}
\end{document}